\documentclass[12pt,a4paper,fleqn]{article}
\usepackage{a4wide,amsfonts,amsmath,latexsym,amssymb,euscript,graphicx,units,mathrsfs}

\usepackage{graphicx}
\usepackage{color}
\usepackage{amssymb}
\usepackage{amssymb}
\usepackage[T1]{fontenc}
\usepackage{latexsym}
\usepackage{xypic}
\usepackage{eufrak}
\usepackage{euscript}
\usepackage{amsfonts,amsmath}
\usepackage{verbatim}
\usepackage{fancyhdr}
\usepackage[english]{babel}
\usepackage{mathrsfs}
\usepackage{units}
\usepackage{hyperref}
\usepackage{tikz}
\usepackage{subfloat}
\usepackage[lofdepth,lotdepth]{subfig}

\newtheorem{prop}{Proposition}[section]

\newtheorem{lemme}[prop]{Lemma}
\newtheorem{rem}[prop]{Remark}

\newtheorem{thm}[prop]{Theorem}
\newtheorem{defi}[prop]{Definition}

\newtheorem{as}[prop]{Assumption}

\renewcommand{\geq}{\geqslant}
\def\leq{\leqslant}

\newcommand{\R}{\mathbb{R}}

\def\1{{\mathbf{1}}}

\def\1{{\mathbf{1}}}
\def\0.5{{\frac{1}{2}}}

\newcommand{\fin}
{ \vspace{-0.6cm}
\begin{flushright}
\mbox{$\Box$}
\end{flushright}
\noindent }

\newcommand{\qed}{\nopagebreak\hspace*{\fill}
{\vrule width6pt height6ptdepth0pt}\par}

\usepackage{fullpage}

\usepackage{hyperref}

\renewcommand{\phi}{\varphi}
\renewcommand{\kappa}{\varkappa}

\newcommand{\fact }{\overline}

\newlength{\querylen}
\usepackage{fancybox}

\begin{document}

\begin{center}
{\Large{\bf The Chen-Stein method for Poisson functionals}}
\normalsize
\\~\\ by Giovanni Peccati\footnote{Universit\'e du Luxembourg. Facult\'e des Sciences, de la Technologie et de la Communication: Unit\'e de Recherche en Math\'ematiques. 6, rue Richard Coudenhove-Kalergi, L-1359 Luxembourg. Email: {\tt giovanni.peccati@gmail.com}} \\ {\it  Universit\'e du Luxembourg}\\~\\
\end{center}

{\small \noindent {\bf Abstract}: We establish a general inequality on the Poisson space, yielding an upper bound for the distance in total variation between the law of a regular random variable with values in the integers and a Poisson distribution. Several applications are provided, in particular: (i) to deduce a set of sufficient conditions implying that a sequence of (suitably shifted) multiple Wiener-It\^o integrals converges in distribution to a Poisson random variable, and (ii) to compute explicit rates of convergence for the Poisson approximation of statistics associated with geometric random graphs with sparse connections (thus refining some findings by Lachi\`eze-Rey and Peccati (2011)). This is the first paper studying Poisson approximations on configuration spaces by combining the Malliavin calculus of variations and the Chen-Stein method. \\



\section{Introduction and motivation}

The aim of this paper is to combine two powerful probabilistic techniques, namely the Chen-Stein method (see e.g. \cite{BHJ, ErSurvey}) and the Malliavin calculus of variations (see e.g. \cite{nuaviv, privaultbook}), in order to study Poisson approximations for functionals of general Poisson random measures. One of our principal achievements is a general inequality on the Poisson space (see Theorem \ref{t:main}), assessing the distance in total variation between the law of a Poisson random variable and the law of a (sufficiently regular) integer-valued Poisson functional. As discussed below, a strong motivation comes from applications in stochastic geometry: in particular, our results allow to generalize, explain and refine some recent findings by Lachi\`eze-Rey and Peccati \cite{LRP1}, dealing with the asymptotic behavior of general random graphs on an Euclidean space. Another remarkable application provided in the paper is a Poisson convergence result for sequences of `perturbed multiple integrals' (see Theorem \ref{t:main2}), extending several existing central limit theorems (CLTs) for multiple Wiener-It\^o integrals (see e.g. \cite{PSTU, PecZheng} and \cite{nuapec}, respectively, for statements in the Poisson and in the Gaussian frameworks). 

Historically, the work \cite{np-ptrf} has been the first paper combining Stein's method for normal approximations (see \cite{ChenGoldShao}) with Malliavin calculus on a Gaussian space. This reference has been the seed of many generalizations and applications, for instance to fractional processes, density estimates and harmonic analysis of Gaussian-subordinated fields on homogenous spaces. See e.g. \cite{np-book} for a presentation of the general theory, and \cite{mp-book} for several applications to the statistical analysis of spherical random fields.

The present work is a natural continuation of the findings contained in \cite{PSTU}, where the authors combined Stein's method with a version of the Malliavin calculus on the Poisson space (due to Nualart and Vives \cite{nuaviv}) in order to compute explicit bounds for CLTs involving functionals of general Poisson measures; see also \cite{PecZheng} for analogous statements in a multi-dimensional setting. The two works \cite{PSTU, PecZheng} have recently triggered many applications in stochastic geometry. The most notable papers in this respect are the following:  reference \cite{lesmathias} lays the foundations of a general asymptotic theory for geometric $U$-statistics; reference \cite{schulte2011} deals with CLTs for Poisson-Voronoi approximations; reference \cite{SchTh} focuses on intrinsic volumes of Poisson $k$-flat processes; reference \cite{DFR} provides an analysis of statistics associated with geometric random graphs on a torus; finally, reference \cite{minh} proves several multidimensional generalizations.

In \cite{LRP1}, the analysis started in \cite{lesmathias} has been extended in order to prove several necessary and sufficient criteria for the normal approximation of random variables having a finite chaotic expansion (based on the use of contraction operators -- see Section \ref{ss:stars} below), as well as to deduce an exhaustive asymptotic characterization of the edge-counting statistics associated with general stationary graphs. In the next subsection, we shall present an overview of the results established in \cite{LRP1} that have provided the impetus for the analysis developed in the present paper.

\begin{rem}{\rm

Other references related to the content of this paper are \cite{KNPS} and  \cite{np-freepoiss}, dealing respectively with the semicircular and the Marchenko-Pastur convergence in the setting of free probability. Note that the Marchenko-Pastur distribution is the free analogous of the Poisson law in classical probability theory. See also \cite{tetilla} for some non-standard free limit theorems, as well as a striking connection with the spanish cheese industry.

}
\end{rem}

\subsection{Motivation: edge counting in random geometric graphs}

Consider an integer $d\geq 1$, and fix the following notation: 
\[
W := \left[-\frac12,\frac12\right]^d, \quad \hat W := \left[-1,1\right]^d, \quad \check{ W }:= \left[-\frac14,\frac14\right]^d
\]
Given two mappings $\lambda \to a(\lambda), \, b(\lambda)$, $\lambda >0$, we write $a(\lambda) \asymp b(\lambda)$ if there exist two positive constants $C,C' >0$ such that $Cb(\lambda) \leq a(\lambda) \leq C' b(\lambda)$ for $\lambda$ sufficiently large. For every $\lambda >0$, we introduce the following objects and notation:

\begin{itemize}

\item[ -- ] $\eta_\lambda$ is a Poisson measure on $(Z,\mathscr{Z} ):=(\R^d, \mathscr{B}(\R^d))$, with control $\lambda \ell$, where $\ell$ is the usual Lebesgue measure on $\R^d$ (see Section \ref{s:framework} for precise definitions).

\item[ -- ] $H_\lambda $ is a symmetric subset of $Z\times Z$ such that: (i) $\ell^2(H_\lambda)<\infty$, and (ii) there exist $\fact H_\lambda \subset Z$ such that $0\notin \fact H_\lambda$ and $H_\lambda = \{(x,y)\in Z\times Z : x-y\in \fact H_\lambda\}$.

\item[ -- ] $F_\lambda = \sum_{x,y\in \eta_\lambda \cap W } {\bf 1}_{H_\lambda} (x,y)$, $F^{\star}_\lambda := \frac12 F_\lambda$ and $\tilde{F}_\lambda := (F_\lambda - E[F_\lambda]) /\sqrt{{\rm Var}(F_\lambda)}$. Note that the random variable $F^{\star}_\lambda$ equals the number of edges in the random graph obtained by taking as vertices the points of the support of $\eta_\lambda$ contained in $W$, and then by connecting two vertices $x,y$ if and only if $(x,y) \in H_\lambda$.  Since $\fact H_\lambda$ does not contain $0$, the graph obtained in this way has no loops, that is, no vertices of the type $\{x,x\}$. Since $H_\lambda$ is symmetric, one has necessarily that $\fact H_\lambda = - \fact H_\lambda$.

\item[ -- ] The three {\it occupation coefficients} introduced in \cite[formulae (4.37)--(4.39)]{LRP1}:
\[
\psi(\lambda) := \ell(\fact H_\lambda\cap W), \quad \hat\psi(\lambda) := \ell(\fact H_\lambda\cap \hat W), \quad \check\psi(\lambda) := \ell(\fact H_\lambda\cap \check W).
\]

\end{itemize}

\begin{rem}{\rm
\begin{enumerate}

\item This class of geometric models contains the so-called {\it Gilbert graphs}, obtained by taking $H_\lambda = \{(x,y)\in \R^d\times \R^d : 0< \|x-y\|_{\R^d} < \delta_\lambda\}$, where $\lambda\to\delta_\lambda$ is a suitable mapping with values in $\R_+$. Graphs corresponding to the cases $d=1$ and $d=2$ are called, respectively, {\it interval graphs} and {\it disk graphs}. See e.g. \cite{CCJ, DFR, penrosebook}.

\item Since $H_\lambda$ can be represented in terms of the set $\fact H_\lambda$, we say that the random graph described above is {\it stationary} -- see again \cite{LRP1}.

\end{enumerate}
}
\end{rem}

The following result is one of the main findings of \cite{LRP1}: it provides an exhaustive description of the asymptotic fluctuations of the family $\{F_\lambda\}$, under an additional regularity assumption on the occupation coefficients $\check\psi, \hat\psi$.

\begin{thm}[See Theorems 4.9-4.11 in \cite{LRP1}] \label{t:lrp}Assume that $\hat\psi(\lambda) \asymp \check\psi(\lambda)$, and consider a standard Gaussian random variable $X\sim \mathscr{N}(0,1)$. The following three properties are in order as $\lambda \to \infty$:
\begin{itemize}

\item[\rm (i)] If $\lambda\psi(\lambda) \to \infty$ or $\lambda\psi(\lambda) \asymp 1$, then there exists a finite constant $K>0$ (independent of $\lambda$) such that $d_W(\tilde F_\lambda,X)\leq K \lambda^{-1/2}$, where $d_W(\tilde F_\lambda,X)$ indicates the Wasserstein distance between the laws of $\tilde F_\lambda$ and $X$, so that $\tilde F_\lambda$ converges to $X$ in distribution. 

\item[\rm (ii)] If $\lambda\psi(\lambda) \to 0$, then ${\rm Var}(F_\lambda) \asymp \lambda^2\psi(\lambda)$. If in addition $\lambda^2\psi(\lambda) \to \infty$, then there exists a constant $K>0$  such that $d_W(\tilde F_\lambda,X)\leq K \left(\lambda \sqrt{\psi(\lambda)}\right)^{-1}$, and consequently $\tilde F_\lambda$ converges to $X$ in distribution. 

\item[\rm (iii)] Assume that $\lambda\psi(\lambda) \to 0$ (so that ${\rm Var}(F_\lambda)\asymp \lambda^2\psi(\lambda)$), and that the mapping $\lambda \to \lambda^2\psi(\lambda)$ is bounded. If there exists a finite constant $c>0$ such that ${\rm Var}(F_\lambda) \to 2c$, then $E[F^\star_\lambda]\to c/2$, and $F^\star_\lambda$ converges in distribution to a Poisson random variable with mean $c/2$ (denoted by ${\rm Po}(c/2)$). 
\end{itemize}

\end{thm}
\begin{rem}{\rm

\begin{enumerate}
\item For explicit examples of sets $H_\lambda$ verifying either one of Points (i)-(iii) in Theorem \ref{t:lrp}, see \cite[Section 4.3.1]{LRP1}. 

\item If the constant $c$ at Point (iii) of the previous statement is equal to zero, then one can prove by a direct argument that $\tilde{F}_\lambda \to 0$ in $L^1$.

\item For every $\lambda >0$, the following equality in law holds:

\begin{equation}
\label{eq:def-F-stat-graphs}
F_{\lambda}\stackrel{\rm Law}{=}\sum_{x,y\in \eta\cap Q_{\lambda}, x\neq y}\1_{x-y\in \fact G_{\lambda}}, \lambda>0,
\end{equation}
where $\eta$ is a random Poisson measure with Lebesgue intensity, and $\fact G_{\lambda}$ is a measurable subset of $\mathbb{R}^d$ defined by the relation 
\begin{equation}\label{e:glambda}
\fact H_{\lambda}=\lambda^{-1/d}\fact G_{\lambda}.
\end{equation}
\end{enumerate}
}
\end{rem}

Point (iii) in Theorem \ref{t:lrp} was proved in \cite{LRP1} by the method of moments, and no information was given about the associated rate of convergence towards the limiting Poisson distribution. In Section \ref{s:graphs}, we will apply the main estimates established in this paper in order to compute an explicit upper bound for the quantity $d_{TV}(F^{\star}_\lambda, {\rm Po}(c/2))$, where $ d_{TV}$ indicates the total variation distance between the laws of $F^{\star}_\lambda$ and ${\rm Po}(c/2)$.

\medskip

The remainder of the paper is organized as follows. Section \ref{s:framework} contains some preliminary results concerning stochastic analysis on the Poisson space. Section \ref{s:main} deals with the main estimates proved in the paper. Section \ref{s:mwii} contains an application to the Poisson approximation of multiple Wiener-It\^o integrals. Finally, in Section \ref{s:graphs} an explicit bound for Point (iii) of Theorem \ref{t:lrp} is computed.

\medskip

\noindent{\bf Acknowledgment.} I am grateful to Rapha\"el Lachi\`eze-Rey for many fundamental conversations about the topics studied in this paper.

\section{Framework}\label{s:framework}

In what follows, the triple $(Z,\mathscr{Z},\mu) $ indicates a measure space such that
$Z$ is a Borel space, $\mathscr{Z}$ is the associated Borel $\sigma$-field, and $\mu$ is a non-atomic $\sigma$-finite
Borel measure. We define $\mathscr{Z}_{\mu} := \{ B\in \mathcal{Z}: \mu(B)< \infty \}$. Throughout the paper, we write
$\eta = \{\eta(B) : B\in \mathscr{Z}_{\mu} \} $ to indicate a {\it Poisson measure} on $(Z,\mathscr{Z}) $ with {\sl control} $\mu$. This means that $\eta $ is a collection of random variables defined on some probability space $(\Omega, \mathscr{F}, P) $, indexed by
the elements of $\mathscr{Z}_{\mu} $ and such that: (a) for every $B,C \in \mathcal{Z}_{\mu}$ such that $B \cap C = \varnothing$, the random variables $ \eta(B)$ and $ \eta(C)$ are independent, and  (b) for every $B \in \mathscr{Z}_{\mu} $, $\eta(B)$ has a Poisson distribution with mean $\mu(B)$. We shall often write $\hat{\eta}(B) = \eta(B) - \mu(B)$, $B\in \mathscr{Z}_\mu$, and $\hat{\eta} = \{\hat{\eta}(B) : B\in \mathscr{Z}_{\mu} \}$. The reader is referred e.g. to \cite{PeTa} for a general introduction to random measures of the Poisson type.  By a slight abuse of notation, we shall sometimes write $ x \in \eta$ in order to indicate that the point $x\in Z$ is charged by the random measure $\eta(\cdot)$.

\begin{rem} \label{rmk1}

\rm{
By virtue of the specific structure of the space $(Z,\mathscr{Z},\mu)$, we can assume throughout the paper that $(\Omega,\mathscr{F},P)$ and $\eta$ are such that
\[ \Omega = \left\{ \omega = \sum_{j=1}^{n} \delta_{z_j},n\in   \mathbb{N} \cup \{\infty\},z_j\in Z    \right\}, \]
where $\delta_z$ denotes the Dirac mass at $z$, and $\eta$ is defined as the
\textit{canonical mapping}:
\[(\omega,B) \mapsto \eta(B)(\omega) = \omega(B) ,\quad B\in \mathscr{Z}_{\mu},\quad \omega\in\Omega.\]  
Finally, the $\sigma$-field $\mathscr{F}$ will be always supposed to be the $P$-completion of the $\sigma$-field generated by $\eta$. }
\end{rem}

\subsection{Chaos and Malliavin calculus}

Given a real $p\in [1,\infty)$, the symbol $L^p(\mu)$ is shorthand for $L^p(Z,\mathscr{Z},\mu)$.
For an integer $q\geq 2$, we shall write $L^p(\mu^q) := L^p(Z^q, \mathscr{Z}^{\otimes q}, \mu^{q}) $, whereas  $L^p_s(\mu^q)$ stands for the subspace of $L^p(\mu^q)$ composed of functions that are $\mu^{q}$-almost everywhere symmetric. We also adopt the convention $L^p(\mu) = L_s^p(\mu) =L^p(\mu^1) =L_s^p(\mu^1) $ and use the following notation: for every $q\geq 1$ and every $f,g\in L^2(\mu^q)$,
  $$ \langle f,g \rangle_{L^2(\mu^q)} = \int_{Z^q} f(z_1,...,z_q)g(z_1,...,z_q)\mu^q (dz_1,...,dz_q), \quad \|f\|_{L^2(\mu^q)} = \langle f,f \rangle^{1/2}_{L^2(\mu^q)} . $$
For every $f\in L^2(\mu^q)$, we denote by $\widetilde{f}$ the {\it canonical symmetrization} of $f$. Plainly, $\|\tilde{f}\|_{L^2(\mu^q)} \leq  \|f\|_{L^2(\mu^q)}$.

\begin{defi}\label{d:mwii}{\rm
For every deterministic function $h\in L^2(\mu)$, we write
\[I_1(h)=\hat{\eta}(h) = \int_Z h(z) \hat{\eta}(dz) \] to indicate the {\it Wiener-It\^o
integral} of $h$ with respect to $\hat{\eta}$. For every $q\geq 2$ and every $f\in L_s^2(\mu^q)$, we denote by $I_q(f)$
the {\it multiple Wiener-It\^o integral}, of order $q$, of $f$ with respect to $\hat{\eta}$. We also set $I_q(f)=I_q(\tilde{f})$, for every $f\in L^2(\mu^q)$ (not necessarily symmetric), and $I_0(b)=b$ for every real constant $b$.
}
\end{defi}

The reader is referred for instance to \cite[Chapter 5]{PeTa} for an exhaustive discussion of multiple Wiener-It\^o integrals and their properties (including the forthcoming Proposition \ref{P : MWIone} and Proposition \ref{P: MWIchaos}).

\begin{prop}\label{P : MWIone}
The following equalities hold for every $q,m\geq 1$, every $f\in L_s^2(\mu^q)$  and every $g\in L_s^2(\mu^m)$:
\begin{enumerate}
  \item[\rm 1.] $E[I_q(f)]=0$,
  \item[\rm 2.] $E[I_q(f) I_m(g)]= q!\langle f,g  \rangle_{L^2(\mu^q)} \1_{(q=m)} $
  {\rm (isometric property).}
\end{enumerate}
\end{prop}

The Hilbert space composed of the random variables with the form $I_q(f)$, where $q\geq 1$ and $f\in L^2_s(\mu^q)$, is called the $q$th \emph{Wiener chaos} associated with the Poisson measure $\eta$. The following well-known {\it chaotic representation property} is an fundamental feature of Poisson random measures. Recall that $\mathscr{F}$ is assumed to be generated by $\eta $.

\begin{prop}
[Wiener-It\^o chaotic decomposition] \label{P: MWIchaos} Every random variable \[F\in L^2(\Omega, \mathscr{F},P):=L^2(P)\]
admits a (unique) chaotic decomposition of the type
\begin{equation} \label{chaos}
F= E[F] + \sum_{i = 1}^{\infty} I_i(f_i),
\end{equation}
where the series converges in $L^2(P)$ and, for each $i\geq 1$, the kernel $f_i$ is an element
of $L^2_s(\mu^i)$.
\end{prop}

\medskip

For the rest of the paper, we shall use definitions and results associated with the Malliavin-type operators defined on the space of functionals of the Poisson measure $\eta$. Our formalism coincides with the one introduced by Nualart and Vives in \cite{nuaviv}. A more recent introduction to the Malliavin calculus on the Poisson space can be found in Privault's monograph \cite{privaultbook}; for the convenience of the reader, some basic definitions and results are presented in the Appendix (see Section \ref{s:appendix}). In particular, we shall denote by $D$, $\delta$, $L$ and $L^{-1}$, respectively, the {\it Malliavin derivative}, the {\it divergence operator}, the {\it Ornstein-Uhlenbeck generator} and its {\it pseudo-inverse}. The domains of $D$, $\delta$ and $L$ are written ${\rm dom} \, D$, ${\rm dom}\, \delta$ and ${\rm dom}\, L$. The domain of $L^{-1}$ is given by the subclass of $L^2(P)$ composed of centered random variables. Given a not necessarily centered $F\in L^2(P)$, one sets by convention $L^{-1} F = L^{-1} (F-E(F))$, so that $LL^{-1} F = F-E(F)$ for every $F\in L^2(P)$. Since the underlying probability space $\Omega$ is assumed to be the collection of discrete measures described in Remark \ref{rmk1}, one can meaningfully define the random variable $\omega\mapsto F_z (\omega) =F(\omega + \delta_z),\, \omega \in \Omega, $  for every given random variable $F$ and every $z\in Z$, where $\delta_z$ is the Dirac mass at $z$. One therefore has the following neat representation of $D$ as a {\it difference operator}:
\begin{lemme}\label{l:diff}
For each $F\in {\rm dom} D$,
$$ D_z F(\omega) = F_z(\omega) - F(\omega) ,\,\, \text{a.s.-} P(d\omega),\,\,  \text{a.e.-} \mu(dz). $$
\end{lemme}
A proof of Lemma \ref{l:diff} is given in \cite{nuaviv}.

\section{A general inequality on the Poisson space}\label{s:main}

The following result is the main finding of the paper. As anticipated in the Introduction, the proof makes use of the so-called {\it Chen-Stein method} for Poisson approximations. A classic reference on the subject is the book by Barbour {\it et al.} \cite{BHJ}; more recent references are the two surveys \cite{CDMsurvey, ErSurvey}. Recall that the {\it total variation distance} between the laws of two random variables $X,Y$ with values in $\mathbb{Z}_+ := \{0,1,2,...\}$ is given by
\begin{equation}\label{e:tv}
d_{TV}(X,Y) = \sup_{A\subset \mathbb{Z}_+} | P(X\in A) - P(Y\in A)| = \frac{1}{2}\sum_{k\geq 0} | P(X=k) - P(Y=k)|.
\end{equation}
Of course, the topology induced by $d_{TV}$ on the class of all probability laws on $\mathbb{Z}_+$ is strictly stronger than the topology induced by convergence in distribution. 

Given a function $f : \mathbb{Z}_+ \to \R$, we denote by $\Delta f$ the {\it forward difference} given by $\Delta f(k) :=  f(k+1) - f(k)$, $k=0,1,2,...\,$; we also use the symbol $\Delta^2 f = \Delta(\Delta f) $. Finally, we write $\|f\|_\infty = \sup_{k\in \mathbb{Z}_+} |f(k)|$.

\begin{thm}[Malliavin bounds for Poisson approximations]\label{t:main} Fix $c >0$, and let ${\rm Po}(c)$ indicate a Poisson random variable with mean $c$. Assume that $F\in L^2(P)$ is an element of ${\rm dom}\,D$ such that $E(F)=c $ and $F$ takes values in $\mathbb{Z}_+$. Then,
\begin{eqnarray}\label{e:mp1}
d_{TV}(F,{\rm Po}(c)) &\leq& \frac{1-e^{-c}}{c} E\left| c - \langle DF, -DL^{-1}F \rangle_{L^2(\mu)}\right| \\
&& +\frac{1-e^{-c}}{c^2} E\left[\int_Z \left| D_zF\,(D_zF-1)\,D_zL^{-1}F \right| \mu(dz)\right]\label{e:mp2}\\
&\leq& \frac{1-e^{-c}}{c} \sqrt{E\left[\left( c - \langle DF, -DL^{-1}F \rangle_{L^2(\mu)}\right)^2\right]}\label{e:mp3} \\
&& +\frac{1-e^{-c}}{c^2} E\left[\int_Z \left| D_zF\,(D_zF-1)\,D_zL^{-1}F \right| \mu(dz)\right].\notag
\end{eqnarray}

\end{thm}

\begin{rem}\label{r:ob}{\rm 
\begin{enumerate}

\item Let $F$ be a centered element of ${\rm dom}\, D$ such that $E(F^2) =1$, and let $Z\sim \mathscr{N}(0,1)$ be a standard Gaussian random variable. In \cite[Theorem 3.1]{PSTU}, it is proved that

\begin{eqnarray}\label{e:mg1}
 d_{W}(F, Z) &\leq& E\left| 1 - \langle DF, -DL^{-1}F \rangle_{L^2(\mu)}\right|\\ 
 &&\quad\quad \quad\quad\quad\quad+ E\left[\int_Z \left| (D_zF)^2\,D_zL^{-1}F \right| \mu(dz)\right],\label{e:mg2}
\end{eqnarray}
where $d_W$ denotes the Wasserstein distance between the laws of two random variables.
\item If $F = \eta(A)$, where $\mu(A) = c$, then one has that $D_z F = -D_zL^{-1}F = {\bf 1}_A(z)$, and therefore \[ c - \langle DF, -DL^{-1}F \rangle_{L^2(\mu)} = \int_Z \left| D_zF\,(D_zF - 1)\,D_zL^{-1}F \right| \mu(dz) = 0,\]
as expected.
\item Let $c, c' >0$. Standard computations (see e.g. \cite[Corollary 3.1]{AL}) yield that
\begin{equation}\label{e:popo} 
d_{TV}({\rm Po}(c), {\rm Po}(c')) \leq |c-c'|.
\end{equation}
It follows that the content of Theorem \ref{t:main} can be extended to a random variable $F$ with arbitrary positive expectation $E(F) = c'$ by using the triangular inequality: \[d_{TV}(F,{\rm Po}(c)) \leq d_{TV}({\rm Po} (c'),{\rm Po}(c))+ d_{TV}(F,{\rm Po}(c'))\leq  |c - c'| +d_{TV}(F,{\rm Po}(c')) .\]
\item Integrating by parts, one sees that $E[\langle DF, -DL^{-1}F \rangle_{L^2(\mu)}] = {\rm Var}(F)$.
\end{enumerate}

}
\end{rem}

\noindent{\it Proof of Theorem \ref{t:main}}. For every $A\subset \mathbb{Z}_+$, we denote by $f_A : \mathbb{Z_+} \to \R$ the unique solution to the Chen-Stein equation
\[
{\bf 1}_A(k) - P({\rm Po}(c) \in A) = c f(k+1) - kf(k), \quad k=0,1,...,
\]
verifying the boundary condition $\Delta^2 f(0) = 0$. Combining e.g. \cite[Theorem 2.3]{ErSurvey} with \cite[Theorem 1.3 and pp. 583-584]{daly}, we immediately deduce the estimates
\begin{equation}\label{e:steinest}
\| f_A\|_\infty \leq \min\left(1, \sqrt{\frac{2}{ce}}\right), \quad \|\Delta f_A\|_\infty \leq \frac{1-e^{-c}}{c}, \mbox{\ and \ } \|\Delta^2 f_A\|_\infty \leq \frac{2-2e^{-c}}{c^2}.
\end{equation}
Using the relation $\delta D = -L$ one infers that
\begin{eqnarray*}
&& P({\rm Po}(c) \in A) - P(F\in A) = E[Ff_A(F) - c f_A(F+1)] \\
&&= E[(F-c)f_A(F) - c \Delta f_A(F)]  = E[\delta(-DL^{-1}F) f_A(F) -  c \Delta f_A(F)].
\end{eqnarray*}
Integrating by parts yields
\[
E[\delta(-DL^{-1}F) f_A(F)] = E[\langle Df_A(F) , -DL^{-1}F \rangle_{L^2(\mu)}],
\]
where, by virtue of Lemma \ref{l:diff}, $D_zf_A(F) = f_A (F+D_zF) - f_A(F)$. Observe that Lemma \ref{l:diff} implies that, since $F$ takes values in $\mathbb{Z}_+$, then one can always choose a version of $D_zF$ with values in $\mathbb{Z}$, in such a way that $F+D_zF = F_z$ takes values in $\mathbb{Z}_+$. Now, for every $f : \mathbb{Z}_+\to \R$ and every $k,a \in \mathbb{Z}_+$ such that $k>a$, one has that
\[
f(k) = f(a) +\Delta f(a) (k-a) +\sum_{j=a}^{k-1} \Delta^2f(j) (k-1-j);
\]
on the other hand, when $k,a \in \mathbb{Z}_+$ are such that $k< a$,
\[
f(k) = f(a) +\Delta f(a) (k-a) +\sum_{j=k}^{a-1} \Delta^2f(j) (j+1-k).
\]
These two relations yield that, for every $k,a \in \mathbb{Z}_+$,
\[
|f(k) - f(a) -\Delta f(a) (k-a)| \leq \frac{\|\Delta^2 f\|_\infty}{2} \left| (k-a)(k-a-1)\right|.
\]
Taking $a =F$ and $k = F_z$, one therefore deduces that
\[
D_z f_A(F) = \Delta f_A(F) D_z F +R_z,
\]
where $R_z$ is a residual random function verifying
\[
|R_z | \leq \frac{\|\Delta^2 f_A\|_\infty}{2} \left| D_z F\,(D_zF-1) \right|, \quad z\in Z.
\]
As a consequence,
\begin{eqnarray*}
&& P({\rm Po}(c) \in A) - P(F\in A) = E[\langle Df_A(F) , -DL^{-1}F \rangle_{L^2(\mu)}- c \Delta f_A(F)]\\
&& = E\big[\Delta f_A(F) (\langle DF , -DL^{-1}F \rangle_{L^2(\mu)}- c )\big ] - E\int_Z ( R_z \times D_zL^{-1}F)\, \mu(dz),
\end{eqnarray*}
and the desired conclusion follows by taking absolute values on both sides, as well as by applying the estimates (\ref{e:steinest}). Inequality (\ref{e:mp3}) follows the Cauchy-Schwarz inequality.
\fin

The following statement is an immediate consequence of Theorem \ref{t:main}. 

\begin{prop}[Poisson limit theorems]\label{p:sufficient}
Fix $c>0$. Let $\{F_n : n\geq 1\}\subset {\rm dom}\, D$ be a sequence of random variables with values in $\mathbb{Z}_+$ such that $E[F_n]\to c$, as $n\to \infty$. Assume that, as $n\to \infty$, 
\begin{enumerate}
\item[\rm 1.] $E\left| c - \langle DF_n, -DL^{-1}F_n \rangle_{L^2(\mu)}\right| \to 0$, and

\item[\rm 2.] $E\left[\int_Z \left| D_zF_n\,(D_zF_n-1)\,D_zL^{-1}F_n \right| \mu(dz)\right] \to 0$.

\end{enumerate}

Then, $\lim_{n\to \infty} d_{TV}(F_n, {\rm Po}(c))\to 0$ and $F_n$ converges in distribution to ${\rm Po}(c)$.

\end{prop}

\begin{rem}{\rm
It is well known that the Poisson distribution is determined by its moments (see e.g. \cite[p. 42-43]{PeTa}). It follows that, in order to prove that a given sequence $\{F_n\}$ converges in distribution to ${\rm Po}(c)$, it is sufficient to prove that $E[F_n^k] \to E[{\rm Po}(c)^k]$, for every integer $k\geq 1$. This is the so called `method of moments', requiring a determination of all moments associated with each $F_n$. Although popular in a Gaussian setting (see \cite{PeTa} for an overview), such a technique is extremely demanding (and very little used) in the framework of Poisson measures. This is mainly due to the fact that the combinatorial structures involved in the so-called `diagram formulae' (that are mnemonic devices used to compute moments by means of combinatorial enumerations -- see \cite[Chapter 4]{PeTa}) become quickly too complex to be effectively put into use. One should compare this situation with the statement of Proposition \ref{p:sufficient}, which only involves two sequences of mathematical expectations. In the forthcoming Sections \ref{s:mwii}--\ref{s:graphs} we will see that, in many applications, in order to check Points 1 and 2 in Proposition \ref{p:sufficient} one is naturally led to assess and control expressions that are only related to the {\it second, third and fourth moments} of the sequence $\{F_n\}$.
}
\end{rem}

\section{Applications to multiple Wiener-It\^o integrals}\label{s:mwii}

In this section, we use Theorem \ref{t:main} in order to establish Poisson convergence results for general sequences of $\mathbb{Z}_+$--valued random variables of the type 
\begin{equation}\label{e:sequence}
F_n = x_n +B_n +I_q(f_n), \quad n\geq 1,
\end{equation}
where: (i) $\{x_n: n\geq 1\}$ is a sequence of positive real numbers, (ii) $q\geq 2$ is an integer independent of $n$, (iii) $I_q$ indicates a multiple Wiener-It\^o integral of order $q$, with respect to the compensated measure $\hat{\eta}$, where $\eta$ is a Poisson measure on the Borel space $(Z,\mathscr{Z})$ with $\sigma$-finite and non-atomic control measure $\mu$, (iv) $f_n \in L_s^2(\mu^q)$, and (v) $\{B_n :n\geq 1\}$ is a {\it smooth vanishing perturbation}, in the sense of the forthcoming Definition \ref{d:svs}.

\begin{rem}{\rm 

\begin{enumerate}
\item In the statement of the forthcoming Theorem \ref{t:main2}, we shall implicitly allow that the underlying Poisson measure $\eta$ also changes with $n$. In particular, one can assume that the associated control measure $\mu = \mu_n$ explicitly depends on $n$. As discussed in Section \ref{s:graphs}, this general framework is useful for geometric applications.

\item We consider perturbed sequences of multiple integrals because, in general, it is not clear whether a non-trivial random variable of the type $x+I_q(f)$, where $x\in \R_+$ and $q\geq 2$, can take values in $\mathbb{Z}_+$, and therefore whether Theorem \ref{t:main} can be directly applied. However, as seen e.g. in Section \ref{s:graphs}, in applications one often encounters sequences of integer-valued random variables whose chaotic decomposition is such that all terms except one vanish asymptotically.
\end{enumerate}

}
\end{rem}

\begin{defi}[Smooth vanishing perturbations]\label{d:svs}{\rm A sequence $\{B_n : n\geq 1\} \subset L^2(P) $ is called a {\it smooth vanishing perturbation} if $B_n, L^{-1}B_n \in {\rm dom}\, D$ for every $n\geq 1$, and the following properties hold:
\begin{eqnarray}
&& \lim_{n\to\infty}E[B_n ^2] = 0\label{e:svs0}\\
\label{e:svs1}
&&\lim_{n\to\infty}E\left[\|DB_n\|^2_{L^2(\mu)}\right] = \lim_{n\to\infty}E\left[\|DL^{-1}B_n\|^2_{L^2(\mu)}\right] =0,\\  
&&\lim_{n\to\infty} E\left[\|DB_n\|^4_{L^4(\mu)}\right] =\lim_{n\to\infty} E\left[\|DL^{-1}B_n\|^4_{L^4(\mu)}\right] =  0.\label{e:svs2}
\end{eqnarray}
Note that, if (\ref{e:svs1})--(\ref{e:svs2}) are verified, an application of the Cauchy-Schwarz inequality yields that
\[
\lim_{n\to\infty}E\left[\|DB_n\|^3_{L^3(\mu)}\right] = \lim_{n\to\infty}E\left[\|DL^{-1}B_n\|^3_{L^3(\mu)}\right] =0
\]
}
\end{defi}

\begin{rem}{\rm Using a Mehler-type representation of the Ornstein-Uhlenbeck semigroup (such as the one stated in \cite[Lemma 6.8.1]{privaultbook}), one sees that the following inequalities are always verified:
\[
  E\left[\|DB_n\|^2_{L^2(\mu)}\right] \geq E\left[\|DL^{-1}B_n\|^2_{L^2(\mu)}\right], \quad   E\left[\|DB_n\|^4_{L^4(\mu)}\right] \geq E\left[\|DL^{-1}B_n\|^4_{L^4(\mu)}\right]. 
\]
}
\end{rem}

\begin{rem}{\rm

The following conditions (a)--(b) are sufficient in order for a given sequence $\{B_n\}$ to be a smooth vanishing perturbation.

\begin{itemize}
\item[\rm (a)] There exists an integer $M\geq 1$, independent of $n$, such that 
\[
B_n = \sum_{i=1}^M I_i(g_{i,n}) :=\sum_{i=1}^M B_{i,n} , \quad n\geq 1,
\]
where $g_{i,n} \in L_s^2(\mu^i)$ for every $i=1,...,M$. Note that these assumptions imply that $E[B_n] = 0$ and $B_n, B_{i,n} \in {\rm dom}\, D$ for every $n$ and every $i$.

\item[\rm (b)] As $n \to \infty$, for every $i=1,...,M$,

\begin{equation}\label{e:svs3}
E[B_n ^2] \to 0, \quad E\left[\int_Z (D_z B_{i,n} )^2\mu(dz)\right] \to 0, \mbox{ and  }  E\left[\int_Z (D_z B_{i,n} )^4\mu(dz)\right] \to 0.
\end{equation}

\end{itemize}
}
\end{rem}

Explicit examples of smooth vanishing perturbations appear in Section \ref{s:graphs} -- see in particular Remark \ref{r:obvious}.

\subsection{Digression: stars and products}\label{ss:stars}

In order to state and prove the main results of this section, it is now necessary to introduce {\it contraction operators} and to discuss their role in product formulae.

The kernel $f \star_r^l g$ on $Z^{p+q-r-l}$, associated with functions $f\in L^2_s(\mu^p) $ and $g \in L^2_s(\mu^q) $, where $p,q \geq 1$, $r=1,\ldots, p\wedge q$ and $l=1,\ldots,r $, is defined as follows:
\begin{eqnarray}
& & f \star_r^l
g(\gamma_1,\ldots,\gamma_{r-l},t_1,,\ldots,t_{p-r},s_1,,\ldots,s_{q-r}) \label{contraction} \\
&=& \int_{Z^l} \mu^l(dz_1,...,dz_l)
f(z_1,,\ldots,z_l,\gamma_1,\ldots,\gamma_{r-l},t_1,,\ldots,t_{p-r}) \nonumber \\
& & \quad\quad\quad\quad\quad\quad\quad\quad\quad\quad\quad\quad \times g(z_1,,\ldots,z_l,\gamma_1,\ldots,\gamma_{r-l},s_1,,\ldots,s_{q-r}). \nonumber
\end{eqnarray}
As it is evident, the operator `$\,\star_r^l\,$' reduces the number of variables in the tensor product of $f$ and $g$ from $p+q$ to $p+q-r-l$: this reduction is obtained by first identifying $r$ variables in $f$ and $g$, and then by integrating out $l$ among them. To deal with the case $l=0$ for $r=0,\ldots, p\wedge q$, we set
\begin{eqnarray*}
& &f \star_r^0
g(\gamma_1,\ldots,\gamma_{r},t_1,,\ldots,t_{p-r},s_1,,\ldots,s_{q-r}) \\
&=& f(\gamma_1,\ldots,\gamma_{r},t_1,,\ldots,t_{p-r})
g(\gamma_1,\ldots,\gamma_{r},s_1,,\ldots,s_{q-r}),
\end{eqnarray*}
 and
$$ f \star_0^0 g (t_1,,\ldots,t_{p},s_1,,\ldots,s_{q}) =f\otimes g (t_1,,\ldots,t_{p},s_1,,\ldots,s_{q})= f(t_1,,\ldots,t_{p})
g(s_1,,\ldots,s_{q}). $$
By using the Cauchy-Schwarz inequality, one sees immediately that $f \star_r^r g$ is square-integrable for any choice of
$r=0,\ldots, p\wedge q$ , and every $f\in L^2_s(\mu^p) $, $g \in L^2_s(\mu^q) $. \\

\begin{rem}\label{r:fub}{\rm
For every
$1\leq p\leq q$ and every
$r=1,...,p$,
\begin{equation}\label{useful}
\int_{Z^{p+q-r}} (f\star_r^0 g)^2 d\mu^{p+q-r} = \int_{Z^r}
(f\star_p^{p-r} f)(g\star_q^{q-r} g) d\mu^r,
\end{equation}
for every $f\in L_s^2(\mu^p)$ and every $g\in L_s^2(\mu^q)$}
\end{rem}

The next result is a fundamental {\it product formula} for Poisson multiple integrals (see e.g. \cite{PeTa} for a proof).
\begin{prop}
[Product formula] Let $f\in L^2_s(\mu^p) $ and $g\in
L^2_s(\mu^q)$, $p,q\geq 1 $, and suppose moreover that $f \star_r^l g
\in L^2(\mu^{p+q-r-l})$ for every $r=1,\ldots,p\wedge q $ and $
l=1,\dots,r$ such that $l\neq r $. Then,
\begin{equation} \label{e:product}
I_p(f)I_q(g) = \sum_{r=0}^{p\wedge q} r!
\left(
\begin{array}{c}
  p\\
  r\\
\end{array}
\right)
 \left(
\begin{array}{c}
  q\\
  r\\
\end{array}
\right)
 \sum_{l=0}^r
 \left(
\begin{array}{c}
  r\\
  l\\
\end{array}
\right)  I_{p+q-r-l} \left(\widetilde{f\star_r^l g}\right),
\end{equation}
 with the tilde $\sim$ indicating a symmetrization.
\end{prop}

In order to be able to directly apply the computations contained in \cite[Proof of Theorem 4.2]{PSTU}, in what follows we shall always work under the following technical assumption.

\begin{as}[Assumption on integrands]\label{a:tech}{\rm Every random variable of the type $Y = I_q(f)$, where $q\geq 2$ and $f\in L^2_s(\mu^q)$, considered in the sequel of this paper is such that the following properties (i)-(iii) are verified.
\begin{enumerate}

\item[(i)] For every $r=1,...q$, the kernel $f_i\star_{q}^{q-r} f_i $ is an element of $L^2(\mu^{r})$.

\item[(ii)] Every contraction of the type $(z_1,...,z_{2q - r- l})\mapsto |f_i|\star_r^l |f_i| (z_1,...,z_{2q - r- l})$ is well-defined and finite for every $r=1,...,q$, every $l=1,...,r$ and every $(z_1,...,z_{2q - r- l})\in Z^{2q-r-l}$.

\item[(iii)] For every $k =  1,..., 2q-2$ and every $(r,l)$ verifying $k = 2q -2-r-l$, 
\[
\int_Z \left[\sqrt{ \int_{Z^k} (f(z,\cdot)\star_r^l f(z,\cdot))^2 \,\,d\mu^k  }\,\,\,\right]\mu(dz)<\infty,
\]
where, for every fixed $z\in Z$, the symbol $f(z,\cdot)$ denotes the mapping $(z_1,...,z_{q-1}) \mapsto f(z,z_1,...,z_{q-1})$.
\end{enumerate}

}
\end{as}

\begin{rem}{\rm According to the discussion contained in \cite{PSTU}, Point (i) in Assumption \ref{a:tech} implies that the following properties (a)-(c) are verified:

\begin{enumerate}

\item[(a)] for every $r=1,...,q$ and every $l=1,...,r$, the contraction $f \star_r^l f$ is a well-defined element of $L^2(\mu^{q_i+q_j-r-l})$;

\item[(b)] for every $r=1,...,q$, $f\star_r^0 f$ is an element of $L^2(\mu^{2q- r})$;

\item[(c)] for every $r=1,...,q$, and every $l=1,...,r\wedge (q-1)$, the kernel $f\star_r^l f $ is a well-defined element of $L^2(\mu^{2q-r-l})$.
\end{enumerate}

In particular, the multiplication formula (\ref{e:product}) implies that every random variable $Y$ verifying Assumption \ref{a:tech} is such that $Y^2 \in L^2(P)$, yielding in turn that $E[Y^4]<\infty$. Analogously, one can also show that, under Assumption \ref{a:tech}, the random variable $\langle DY, -DL^{-1}Y\rangle_{L^2(\mu)}$ is square-integrable (and not merely an element of $L^1(P)$).
}
\end{rem}

\begin{rem}{\rm For instance, Assumption \ref{a:tech} is verified whenever $f$ is a bounded function with support in a rectangle of the type $B\times\cdots\times B$, where $\mu(B)<\infty$. 
}
\end{rem}

\subsection{Poisson limit theorems}
 
The following statement contains the main result of the section. 
  
\begin{thm}[Poisson limit theorems on a perturbed chaos]\label{t:main2} Fix $c>0$ and let ${\rm Po}(c)$ be a Poisson random variable with mean $c$.
For a fixed $q\geq 2$, let $\{ F_n : n\geq 1 \}$ be a $\mathbb{Z}_+$--valued sequence as in (\ref{e:sequence}), such that $x_n \to c$ and $E[I_q(f_n)^2]\to c$.  Assume moreover that the following three conditions hold:
\begin{itemize}
\item[\rm (i)] For every $n\geq 1$, the kernel $f_n$ verifies
Assumption \ref{a:tech}.
\item[\rm (ii)] For every $r=1,...,q$, and every $l=1,...,r \!
\wedge \! (q-1)$, $\|f_n \star_r^l f_n
\|_{L^2(\mu^{2q-r-l})} \rightarrow 0$ (as $n\rightarrow \infty$).
\item[\rm (iii)] The relation $\sup_n\|f_n\|_{L^4(\mu^q)}<\infty$ holds and, as  $n\rightarrow\infty$,
\begin{equation}\label{e:myprecious}
\int_{Z^q} \big(\, f_n^2 + q!^2 f^4_n -2q! f_n^3\, \big)\,\, d\mu ^q \longrightarrow 0.
\end{equation}
\end{itemize}
Then, $d_{TV}(F_n, {\rm Po}(c)) \to 0$, as $n\rightarrow
\infty$.

\end{thm}

\begin{rem}\label{r:myprecious}{\rm
\begin{enumerate}

\item Condition (\ref{e:myprecious}) is trivially verified whenever \[f_n(z_1,...,z_q) =\frac{1}{q!}\, {\bf 1}_{H_n}(z_1,...,z_q),\] where $H_n$ is some measurable symmetric subset of $Z^q$.

\item When $q = 2$, Conditions (ii) and (iii) in the statement of Theorem \ref{t:main2} boil down to the following asymptotic relations
\begin{eqnarray*}
&& \|f_n \star_1^1 f_n
\|_{L^2(\mu^{2})} \rightarrow 0, \quad  \|f_n \star_2^1 f_n
\|_{L^2(\mu)} \rightarrow 0, \quad \mbox{and}\\
&& \int_{Z^2} \big(\, f_n^2 + 4 f^4_n -4 f_n^3\, \big)\,\, d\mu ^2 \longrightarrow 0.
\end{eqnarray*} 
\end{enumerate}
}
\end{rem}

Theorem \ref{t:main2} should be compared with the following central limit result, first proved in \cite[Theorem 4.2]{PSTU}.

\begin{thm}[CLTs on a fixed chaos, see \cite{PSTU}]\label{t:chaosCLT} Let $X\sim \mathscr{N}(0,1)$.
Fix $q\geq 2$, and let $F_n =
I_q(f_n)$, $n\geq 1$, be a sequence of multiple stochastic
Wiener-It\^o integrals of order $q$. Suppose that, as
$n\rightarrow\infty$,
$E(F_n^2) \rightarrow 1$. Assume in addition that the following conditions are verified:
\begin{itemize}
\item[\rm (i)] For every $n\geq 1$, the kernel $f_n$ verifies
Assumption \ref{a:tech}.
\item[\rm (ii)] For every $r=1,...,q$, and every $l=1,...,r \!
\wedge \! (q-1)$, $\|f_n \star_r^l f_n
\|_{L^2(\mu^{2q-r-l})} \rightarrow 0$ (as $n\rightarrow \infty$).
\item[\rm (iii)] As $n\rightarrow\infty$,
\begin{equation}\label{e:barcelona}
\int_{Z^q}  f^4_n \,\, d\mu ^q \longrightarrow 0.
\end{equation}
\end{itemize}
Then, $F_n $ converges in distribution to $X$, as $n\rightarrow
\infty$, in the sense of the Wasserstein distance.
\end{thm}

\medskip

\noindent{\it Proof of Theorem \ref{t:main2}}. We have to prove that, if (i)-(iii) are verified, then Point 1 and Point 2 in the statement of Proposition \ref{p:sufficient} hold. The proof is divided into three steps.

\smallskip

\noindent{\it Step 1: Computations related to $DI_q(f_n)$}. One has that $D_zI_q(f_n) = qI_{q-1}(f_n(z,\cdot)$ and \[\langle DI_q(f_n) , -DL^{-1}I_q(f_n)\rangle_{L^2(\mu)} = \frac1q\|DI_q(f_n)\|_{L^2(\mu)}.\] Using the computations contained in \cite[p. 464]{PSTU}, one sees that
\begin{equation}\label{e:dev}
\{D_zI_q(f_n)^2\} =q^2 \sum_{p=0}^{2q-2} I_p(G_p^{q-1}f(z,\cdot)),
\end{equation}
where 
\begin{equation}\label{e:der2}
G_{p}^{q-1}f(z,\cdot)(z_1,...,z_p) = \sum_{r=0}^{q-1}\sum_{l=0}^r {\bf 1}_{\{2q-2-r-l=p\}} r!\binom{q-1}{r}^2\binom{r}{l} \widetilde{f(z,\cdot)\star_r^lf(z,\cdot)}(z_1,...,z_p),
\end{equation}
and the stochastic integrals are set equal to zero on the exceptional set composed of those $z$ such that $f(z,\cdot)\star_r^lf(z,\cdot)$ is not an element of $L^2(\mu^{2q-2-r-l})$ for some $r,l$. Combining the triangular and Cauchy-Schwarz inequalities with the estimates in \cite[Theorem 4.2]{PSTU} (see also \cite[Theorem 3.5]{LRP1} for a more compact statement), one deduces that there exists a constant $K$, depending uniquely on $c$ and $q$, such that, writing $E[I_q(f_n)^2]:= y_n$,
\begin{eqnarray}\label{e:rough1}
&&E[|c - q^{-1}\|DI_q(f_n)\|_{L^2(\mu)}^2|] \leq |c-y_n| + K\times \max_{\substack{r=1,...,q \\ l=1,...,r\wedge (q-1)}} \|f_n\star_r^l f_n\|_{L^2(\mu^{2q-r-l})}\\
&& \sqrt{E[ \|DI_q(f_n)\|^4_{L^4(\mu)}]} \leq K\left\{\max_{\substack{r=1,...,q \\ l=1,...,r\wedge (q-1)}} \|f_n\star_r^l f_n\|_{L^2(\mu^{2q-r-l})}+ \|f_n\|^2_{L^4(\mu)}\right\}. \label{e:rough2}
\end{eqnarray}
In particular, these relations imply that the sequence $n\mapsto E[\|DI_q(f_n)\|_{L^k(\mu)}^k]$ is bounded for $k=2,3,4$.
\smallskip

\noindent{\it Step 2: Dealing with $B_n$}. Since $\{B_n\}$  is a smooth vanishing perturbation, using the Cauchy-Schwarz inequality and the results from Step 1, one deduces immediately that the sequence
\begin{eqnarray*}
&&\langle DF_n , -DL^{-1}F_n \rangle_{L^2(\mu)} - \langle DI_q(f_n) , -DL^{-1}I_q(f_n) \rangle_{L^2(\mu)}\\
&& =\langle DF_n , -DL^{-1}F_n \rangle_{L^2(\mu)} -\frac1q \| DI_q(f_n) \|^2_{L^2(\mu)}  \\
&& = \langle DI_q(f_n) , -DL^{-1}B_n \rangle_{L^2(\mu)}+\langle DB_n , -DL^{-1}I_q(f_n) \rangle_{L^2(\mu)}+ \langle DB_n , -DL^{-1}B_n \rangle_{L^2(\mu)}
\end{eqnarray*}
converges to zero in $L^1(P)$. Exploiting in a similar way the fact that $\{B_n\}$ is a smooth vanishing perturbation together with the bounded character of $E[\|DI_q(f_n)\|_{L^2(\mu)}^2]$ and $E[\|DI_q(f_n)\|_{L^4(\mu)}^4]$, one also deduces that 
\begin{eqnarray*}
&& E\left[\left| \int_Z \left| D_zF_n\,(D_zF_n-1)\,D_zL^{-1}F_n \right| \mu(dz)\right.\right. \\
&& \quad\quad\quad \left.\left.- \int_Z \left| D_zI_q(f_n)\,(D_zI_q(f_n)-1)\,D_zL^{-1}I_q(f_n) \right| \mu(dz)\right|\right] \to 0.
\end{eqnarray*}
Now observe that, thanks to (\ref{e:rough1}) and Assumption (ii) in the statement,
\[
E[|c - q^{-1}\|DI_q(f_n)\|_{L^2(\mu)}^2|] \to 0, \quad n\to\infty.
\] 
It follows that the proof is concluded once we prove that, as $n \to \infty$,
\begin{eqnarray}
&&E\left[ \int_Z \left| D_zI_q(f_n)\,(D_zI_q(f_n)-1)\,D_zL^{-1}I_q(f_n) \right| \mu(dz)\right]\notag \\
&&\quad\quad\quad =q^{-1}E\left[ \int_Z  (D_zI_q(f_n))^2\,|D_zI_q(f_n)-1|  \mu(dz)\right] \to 0.\label{e:step3}
\end{eqnarray}
This is the object of the forthcoming Step 3.
\smallskip

\noindent{\it Step 3: Proving (\ref{e:step3})}. Using the Cauchy-Schwarz inequality, we infer that
\begin{eqnarray*}
&& E\left[ \int_Z  (D_zI_q(f_n))^2\,|D_zI_q(f_n)-1|  \mu(dz)\right]\\
&& \leq \sqrt{E[\|DI_q(f_n)\|_{L^4(\mu)}^4]}\times \sqrt{ \int_Z  E\left[(D_zI_q(f_n))^2\,(D_zI_q(f_n)-1)^2 \right] \mu(dz)}
\end{eqnarray*}
Also, recall that $q^2G_0^{q-1}f_n(z,\cdot) = qq!\int_{Z^{q-1}} f(z,\cdot)^2d\mu^{q-1}$, and consequently
\[
E[\|DI_q(f_n)\|^2_{L^2(\mu)} ]= q q!\|f_n\|^2_{L^2(\mu^q)} .
\]
 From (\ref{e:dev}), we infer that
\begin{eqnarray*}
&& E\left[(D_zI_q(f_n))^2\,(D_zI_q(f_n)-1)^2 \right] = (1+q^2G_0^{q-1}f_n(z,\cdot))\times q^2G_0^{q-1}f_n(z,\cdot) \\
&& + q^4 \sum_{\substack{p=1,...,2q-2 \\ p\neq q-1}}^{2q-2}p! \|\widetilde{G_p^{q-1}f(z,\cdot))}\|^2_{L^2(\mu^p)}\\
&&+ q^4(q-1)! \int_Z \widetilde{G_{q-1}^{q-1}f_n(z,\cdot)}\left(\widetilde{G_{q-1}^{q-1}f_n(z,\cdot)}-2q^{-1}f(z,\cdot)\right) d\mu^{q-1}.
\end{eqnarray*}
Integrating and simplifying the RHS of the previous equality by reasoning as in \cite[p. 467]{PSTU}, one deduces from Assumption (ii) in the statement that
\[
\int_Z E\left[(D_zI_q(f_n))^2\,(D_zI_q(f_n)-1)^2 \right]\mu(dz) = o(1) + \int_{Z^q} \left( qq!f_n^2 + qq!^3 f_n^4 - 2qq!^2 f_n^3 \right) d\mu^q, 
\]
where $o(1)$ indicates a sequence converging to 0, as $n\to\infty$. The conclusion follows by applying relation (\ref{e:myprecious}).
\fin

\section{An application to geometric random graphs}\label{s:graphs}

The following statement provides the announced refinement of Theorem \ref{t:lrp}.

\begin{thm}\label{t:main3} Let the notation and assumptions of Theorem \ref{t:lrp}-(iii) prevail. Then, there exists a finite constant $K$, independent of $\lambda$, such that
\begin{equation}\label{e:gb} 
d_{TV}(F_\lambda^\star, {\rm Po}(c/2))\leq |E[F^\star_\lambda] - c/2| + K\sqrt{\lambda\psi(\lambda)}.
\end{equation}
\end{thm}
{\it Proof.} Using e.g. \cite[Theorem 3.1]{lesmathias}, one sees that $F^\star_\lambda$ admits the following chaotic decomposition
\[
F^\star_\lambda = E[F^\star_\lambda] +I_1(f_{1,\lambda}) +I_2(f_{2,\lambda}),
\]
where $I_1, I_2$ indicate (multiple) Wiener-It\^o integrals with respect to the compensated measure $\hat\eta_\lambda = \eta - \lambda\ell$, $f_{1,\lambda} (z) = \lambda \int_W {\bf 1}_{H_\lambda}(z,x)\ell(dx) {\bf 1}_W(z)$, and $f_{2,\lambda}(z_1,z_2) = \frac12{\bf 1}_{H_\lambda\cap(W\times W)}(z_1,z_2)$. Also, Campbell's Theorem (see \cite[Theorem 3.1.3]{SW}) implies that $E[F^\star_\lambda] =\frac12 \lambda^2 \ell(H_\lambda)$. Applying Theorem \ref{t:main} together with Remark \ref{r:ob} and the Cauchy-Schwarz inequality yields that
\begin{eqnarray*}
&& d_{TV}(F^\star_\lambda, {\rm Po}(c/2)) \leq  |E[F^\star_\lambda] - c/2| +\frac{2-2e^{-c/2}}{c}\Xi_0(\lambda) +\frac{4-4e^{-c/2}}{c^2}\Xi_1(\lambda)\times \Xi_2(\lambda), 
\end{eqnarray*}

where

\begin{eqnarray*}
&& \Xi_0(\lambda):= \sqrt{E\left[\left( E[F^\star_\lambda] - \langle DF^\star_\lambda, -DL^{-1}F^\star_\lambda \rangle_{L^2(\lambda\ell)}\right)^2\right]} \\
&& \Xi_1(\lambda) := \sqrt{\lambda E\left[ \int_Z (D_zL^{-1} F^\star_\lambda)^2\ell(dz)\right]},\\
&& \Xi_2(\lambda) := \sqrt{\lambda E\left[ \int_Z (D_z F^\star_\lambda)^2(D_z F^\star_\lambda-1)^2\ell(dz)\right]}. 
\end{eqnarray*}
One has that $D_z F^\star_\lambda = f_{1,\lambda}(z) + 2I_1(f_{2,\lambda}(z,\cdot))$. Using the computations contained in \cite[p. 470]{PSTU}, we infer that
\begin{eqnarray*}
&&\langle DF^\star_\lambda, -DL^{-1}F^\star_\lambda \rangle_{L^2(\lambda\ell)} \\
&& = \|f_{1,\lambda}\|^2_{L^2(\lambda\ell)}+2\|f_{2,\lambda}\|^2_{L^2((\lambda\ell)^2)} + 2I_1(f_{2,\lambda}\star_2^1 f_{2,\lambda}) + 2I_2(f_{2,\lambda}\star_1^1 f_{2,\lambda})+ 3I_1(f_{1,\lambda}\star_2^1 f_{2,\lambda}).  
\end{eqnarray*}
Using the computations contained in \cite[Theorem 4.8]{LRP1} one sees that the following five relations are verified: 
\begin{enumerate}

\item[(1)]  $\|f_{1,\lambda}\|^2_{L^2(\lambda\ell)}\asymp \lambda^3\psi(\lambda)^2\asymp \lambda\psi(\lambda)$, 
\item[(2)] $2\|f_{2,\lambda}\|^2_{L^2((\lambda\ell)^2) }=E[F_\lambda^\star] $, 
\item[(3)] $\|f_{2,\lambda}\star_2^1f_{2,\lambda}\|^2_{L^2(\lambda\ell)}\asymp \lambda^3\psi(\lambda)^2\asymp \lambda\psi(\lambda)$, 
\item[(4)] $\|f_{2,\lambda}\star_1^1f_{2,\lambda}\|^2_{L^2((\lambda\ell)^2)}\asymp \lambda^4\psi(\lambda)^3\asymp (\lambda\psi(\lambda))^2$,
\item[(5)] $\|f_{1,\lambda}\star_1^1f_{2,\lambda}\|^2_{L^2(\lambda\ell)}\asymp \lambda^5\psi(\lambda)^4\asymp (\lambda\psi(\lambda))^3$.  
\end{enumerate}
This implies in particular that $\Xi_0(\lambda) \asymp \sqrt{\lambda\psi(\lambda)}$. By using the explicit expression $-D_zL^{-1} F^\star_\lambda = f_{1,\lambda}(z) + I_1(f_{2,\lambda}(z,\cdot))$, it is not difficult to prove that the mapping $\lambda\mapsto \Xi_1(\lambda)$ is necessarily bounded, so that the statement is proved once we show that $\Xi_2(\lambda) \asymp \sqrt{\lambda\psi(\lambda)}$. Using the relation 
\begin{eqnarray*}
(D_zF^\star_\lambda) ^2 &=& \left[f_{1,\lambda}(z)^2 +4\lambda\int_Z f_{2,\lambda}(z,x)\ell(dx)\right] \\ 
&& +4I_1[f_{2,\lambda}(z,\cdot)(f_{1,\lambda}(z) + f_{2,\lambda}(z,\cdot))] + 4I_2(f_{2,\lambda}\star_0^0 f_{2,\lambda}),
\end{eqnarray*}
developing the square $(D_zF^\star_\lambda - 1)^2$ and integrating with respect to $z$ yields
\begin{eqnarray*}
&& \lambda E\left[ \int_Z (D_z F^\star_\lambda)^2(D_z F^\star_\lambda-1)^2\ell(dz)\right]\\
&&= \lambda\int_Z  (f_{1,\lambda}(z)^2 +f_{1,\lambda}(z)^4)\ell(dz)+8\lambda^2 \int_Z\int_Z f_{1,\lambda}(z)f_{2,\lambda}(z,x)^2 \ell(dx)\ell(dz)\\
&& +48\lambda^3\int_Z \left(\int_Z f_{2,\lambda}(z,x)^2\ell(dx)\right)^2\ell(dz)+ 4\lambda^2\int_Z\int_Z f_{2,\lambda}(z,x)^2\ell(dx)\ell(dz)\\
&& + 16\lambda^2\int_Z\int_Z  f^2_{2,\lambda}(z,x)(f_{1,\lambda}(z) +f_{2,\lambda}(z,x))(f_{1,\lambda}(z) +f_{2,\lambda}(z,x)-1)\ell(dx)\ell(dz).
\end{eqnarray*}
Exploiting the explicit form of $f_{2,\lambda}$, one sees that
\[
4\lambda^2\int_Z f_{2,\lambda}(z,x)^2\ell(dx)\ell(dz) + 16\lambda^2\int_Z f_{2,\lambda}(z,x)^4\ell(dx)\ell(dz) - 16\lambda^2\int_Z f_{2,\lambda}(z,x)^3\ell(dx)\ell(dz) = 0.
\]
One can now deal with the remaining terms by using once again the estimates contained in \cite[Theorem 4.8]{LRP1}: this entails the relation $\lambda E\left[ \int_Z (D_z F^\star_\lambda)^2(D_z F^\star_\lambda-1)^2\ell(dz)\right]\asymp \lambda\psi(\lambda)$, and consequently the desired conclusion $\Xi_2(\lambda) \asymp \sqrt{\lambda\psi(\lambda)}$.
\fin

\begin{rem}\label{r:obvious}{\rm
By inspection of the previous proof, one sees that the mapping $\lambda \mapsto I_1(f_{1,\lambda})$ is a smooth vanishing perturbation, in the sense of Definition \ref{d:svs}. It follows that the convergence $F_\lambda^\star \stackrel{\rm Law}{\rightarrow} {\rm Po}(c/2)$ can alternatively be proved by directly applying Theorem \ref{t:main2}.
}
\end{rem}

\bibliographystyle{plain}

\begin{section}{Appendix: Malliavin operators on the Poisson space}\label{s:appendix}
We now define some Malliavin-type operators associated with a Poisson measure $\eta$, on the Borel space $(Z,\mathscr{Z})$, with non-atomic control measure $\mu$. We follow the work by Nualart and Vives \cite{nuaviv}.

\medskip

\noindent \underline{\bf The derivative operator $D$}. 

\smallskip

For every $F\in L^2(P)$, the derivative of $F$, $DF$ is defined as an element of $L^2(P;L^2(\mu))$, that is, of the space of the jointly measurable random functions $u:\Omega \times Z \mapsto \mathbb{R}$ such that $E \left[\int_Z u_z^2 \mu(dz) \right] <\infty$.
\begin{defi}
 \begin{enumerate}
   \item The domain of the derivative operator $D$, written  ${\rm dom} D$, is the set of all random variables $F\in L^2(P)$ admitting a chaotic decomposition (\ref{chao}) such that
$$ \sum_{k\geq 1} k k!\|f_k \|^2_{L^2(\mu^k)} < \infty ,$$
   \item For any $F\in {\rm dom}D$, the random function $z \mapsto D_z F$ is defined by
$$ D_z F= \sum_{k \geq 1}^{\infty} k I_{k-1}(f_k(z,\cdot)) .$$
 \end{enumerate}
\end{defi}

\medskip

\noindent \underline{\bf The divergence operator $\delta$}.  

\smallskip

Thanks to the chaotic representation property of $\eta$, every random function
$u \in L^2(P,L^2(\mu))$ admits a unique representation of the type
\begin{equation} \label{skor}
 u_z = \sum_{k \geq 0}^{\infty}  I_{k}(f_k(z,\cdot)) ,\,\, z\in Z,
\end{equation}
where the kernel $f_k$ is a function of $k+1$ variables, and $f_k(z,\cdot)$ is an element of $L^2_s(\mu^k)$. The {\sl divergence operator} $\delta(u)$ maps a random function $u$ in its domain to an element of $L^2(P)$.\\

\begin{defi}
\begin{enumerate}
  \item  The domain of the divergence operator, denoted by  ${\rm dom} \delta$, is the collection of all $u\in L^2(P,L^2(\mu))$ having the above chaotic expansion (\ref{skor}) satisfied the condition:
$$ \sum_{k\geq 0}  (k+1)! \|f_k \|^2_{L^2(\mu^(k+1))} < \infty. $$
  \item For $u\in {\rm dom}\delta$, the random variable $\delta(u)$ is given by
      $$ \delta (u) = \sum_{k\geq 0} I_{k+1}(\tilde{f}_k), $$
      where $\tilde{f}_k$ is the canonical symmetrization of the $k+1$ variables function $f_k$.
\end{enumerate}
\end{defi}
As made clear in the following statement, the operator $\delta$ is indeed the adjoint operator of $D$.
\begin{lemme}[Integration by parts]\label{L : IBP}
 For every $G\in {\rm dom} D$ and $u\in {\rm dom} \delta$, one has that
$$ E[G \delta(u)] = E[\langle D G, u \rangle_{L^2(\mu)}]. $$
\end{lemme}
The proof of Lemma \ref{L : IBP} is detailed e.g. in \cite{nuaviv}.\\

\medskip

\noindent \underline{\bf The Ornstein-Uhlenbeck generator $L$}.

\smallskip

\begin{defi}
\begin{enumerate}
  \item  The domain of the Ornstein-Uhlenbeck generator, denoted by  ${\rm dom} L$, is the collection of all $F \in L^2(P)$ whose chaotic representation \label{chao} verifies the condition:
$$ \sum_{k\geq 1}  k^2 k! \|f_k \|^2_{L^2(\mu^k)} < \infty $$
  \item The Ornstein-Uhlenbeck generator $L$ acts on random variable $F\in {\rm dom}L$ as follows:
      $$ LF = - \sum_{k\geq 1} k I_{k}(f_k) .$$
\end{enumerate}
\end{defi}

\medskip

\noindent \underline{\bf The pseudo-inverse of $L$}.

\smallskip

\begin{defi}
\begin{enumerate}
  \item  The domain of the pseudo-inverse of the Ornstein-Uhlenbeck generator, denoted by $L^{-1}$, is the space $L^2_0(P)$ of \it{centered} random variables in $L^2(P)$.
  \item For $F = \sum\limits_{k\geq 1} I_k (f_k) \in L^2_0(P)$ , we set
      $$ L^{-1}F = - \sum_{k\geq 1} \cfrac{1}{k} I_{k}(f_k). $$
\end{enumerate}
\end{defi}
\end{section}

\end{document}